\documentclass[11pt,a4paper]{article}
\usepackage{mathrsfs}
\usepackage{amsfonts}
\usepackage{amssymb}
\usepackage{color,xcolor}
\usepackage{manfnt}
\usepackage{graphicx}
\usepackage{amsfonts,amsmath, amssymb}

\newtheorem{theo}{\textbf{\ \ \quad Theorem}}[section]

\newtheorem{lem}{\textbf{\ \ \quad Lemma}}[section]
\newtheorem{remark}{\textbf{\ \ \quad Remark}}[section]
\newtheorem{col}{\textbf{\ \ \quad Corollary}}[section]
\newtheorem{prop}{\textbf{\ \ \quad Proposition}}[section]
\newtheorem{defi}{\textbf{\ \ \quad Definition}}[section]

\newcommand{\lbl}[1]{\label{#1}}

\newcommand{\be}{\begin{equation}}
\newcommand{\ee}{\end{equation}}
\newcommand\bes{\begin{eqnarray}}
\newcommand\ees{\end{eqnarray}}
\newcommand{\bess}{\begin{eqnarray*}}
\newcommand{\eess}{\end{eqnarray*}}

\newcommand{\nm}{\nonumber}

\newcommand{\ds}{\displaystyle}

\newcommand{\R}{\mathbb{R}}

{\setlength\arraycolsep{2pt}}

\title{Morrey-Campanato estimates for the moments of stochastic integral operators and its application to SPDEs }

\author{Guangying Lv$^{a,b}$, Hongjun Gao$^b$, Jinlong Wei$^c$, Jiang-Lun Wu$^d$\\
\\
\ \\
   {\small \it $^a$Institute of Contemporary Mathematics, Henan University}\\
  {\small \it Kaifeng, Henan 475001, China}\\
  {\small \tt gylvmaths@henu.edu.cn}\\
    {\small \it $^b$Institute of Mathematics, School of Mathematical Science}\\
  {\small \it Nanjing Normal University, Nanjing 210023, China}\\
  {\small \tt gaohj@njnu.edu.cn}\\
  {\small \it $^c$ School of Statistics and Mathematics, Zhongnan University of}\\
  {\small \it
 Economics and Law, Wuhan, Hubei 430073, China}\\
   {\small \tt  weijinlong.hust@gmail.com }\\
   {\small \it $^d$ Department of Mathematics, Swansea University, Swansea SA2 8PP, UK}\\
   {\small \tt  j.l.wu@swansea.ac.uk }
}

\begin{document}
\maketitle

\medskip

\begin{abstract}
In this paper, we are concerned with the estimates for the moments of stochastic convolution integrals. We first deal with the stochastic singular integral operators and we aim to derive the Morrey-Campanato estimates for the $p$-moments (for $p\ge1$). Then, by utilising the embedding theory between the Campanato space and H\"older space, we establish the norm of $C^{\theta,\theta/2}(\bar D)$, where $\theta\ge0, \bar D=\bar G\times[0,T]$ for arbitrarily fixed $T\in(0,\infty)$ and $G\subset\mathbb{R}^d$. As an application, we consider the following stochastic (fractional) heat equations with additive noises
   \bess
du_t(x)=\Delta^\alpha u_t(x)dt+g(t,x)d\eta_t,\ \ \ u_0=0,\ 0\leq t\leq T, x\in G,
   \eess
where $\Delta^\alpha=-(-\Delta)^\alpha$ with $0<\alpha\leq1$ (the fractional Laplacian),
$g:[0,T]\times G\times\Omega\to\mathbb{R}$ is a joint measurable coefficient, and $\eta_t, t\in[0,T]$, is either the Brownian motion or a L\'evy process on a given filtered probability space $(\Omega,\mathcal{F},P;\{\mathcal{F}_t\}_{t\in[0,T]})$. The Schauder estimate for the $p$-moments of the solution of the above equation is obtained. The novelty of the present paper is that we obtain the Schauder estimate for parabolic stochastic partial differential equations with L\'evy noise.



\end{abstract}

\baselineskip=15pt

\section{Introduction}
\setcounter{equation}{0}
For a stochastic process $\{X_t,t\in [0,T]\}$, there are two important aspects worth investigating.
One is the associated probability density functions (PDFs) or its probability laws, and the other
is the moments estimate. But for a stochastic process depending on spatial variable (to be more precise,
a random field), that is, $X_t=X(t,x,\omega)$ with $x$ standing for a spatial variable, it would be hard to
consider its PDFs or probability laws. Fortunately, we could get some moments estimate. In this paper,
we focus on the estimates for solutions of (parabolic) stochastic partial differential equations (SPDEs),
in particular, on the Schauder estimate for the SPDEs.

For (parabolic) SPDEs, certain kinds of estimates for the solutions have been well studied. By using parabolic
Littlewood-Paley inequality, Krylov \cite{Krylov} proved that for SPDEs of the following type
   \bes
du=\Delta udt+gdw_t,
  \lbl{1.1}\ees
it holds that
   \bes
\mathbb{E}\|\nabla u\|_{L^p((0,T)\times\mathbb{R}^d)}^p\leq C(d,p)
\mathbb{E}\|g\|_{L^p((0,T)\times\mathbb{R}^d)}^p,
    \lbl{1.2} \ees
where $w_t$ is a Wiener process and $p\in[2,\infty)$. Moreover, van Neerven et al. \cite{NVW} made
a significant extension of (\ref{1.2}) to a class of operators $A$ which admit a bounded
$H^\infty$-calculus of angle less than $\pi/2$. Kim \cite{Kim} established a BMO estimate
 for stochastic singular integral operators. And as an application, he considered (\ref{1.1}) and interestingly
 he obtained the $q$-th order BMO quasi-norm of the
derivative of $u$ is controlled by $\|g\|_{L^\infty}$. More recently, Kim et al. \cite{KKL}
studied the parabolic Littlewood-Paley inequality for a class of time-dependent pseudo-differential operators
of arbitrary order, and applied their result to a high-order stochastic PDE. We refer the interested readers to
\cite{Cbook,Gbook} for a comprehensive account on the BMO estimates.

Recently, Yang \cite{Y} considered the following SPDEs
 \bess
du=\Delta^{\frac{\alpha}{2}}udt+fdX_t,\ \ \ u_0=0,\ 0< t<T,
  \eess
where $\Delta^{\frac{\alpha}{2}}=-(-\Delta)^{\frac{\alpha}{2}}$, for $0<\alpha<2$, are nothing but the fractional Laplace operators,
and $X_t$ is a L\'{e}vy process. The author obtained a parabolic Triebel-Lizorkin space estimate for the convolution
operator.

In our paper \cite{LGWW}, we consider the stochastic singular operator
    \bes
\mathcal {G}g(t,x)&=&\int_0^t\int_ZK(t,s,\cdot)\ast g(s,\cdot,z)(x)\tilde N(dz,ds)\nm\\
&=&\int_0^{t}\int_Z\int_{\mathbb{R}^d}K(t-s,x-y)g(s,y,z)dy\tilde N(dz,ds),
   \lbl{1.3}\ees
for $g:[0,T]\times\mathbb{R}^d\times Z\times\Omega\to\mathbb{R}$ being a predictable process,
where $\tilde N$ is a compensated Poisson measure.
Under appropriate conditions on the kernel $K$, we obtained the $q$-th order BMO
estimate. As an application, we obtained the $q$-th order BMO estimate for the solution of
the stochastic nonlocal heat equation.

For the regularity of SPDEs, several important works have been established, see
\cite{K2004,K2008,K1996,NVW,Z2006}. Similar to the regularity of PDE, the regularity of SPDEs can be
divided into two parts. One is the $L^p$-theory. Krylov \cite{K1996} obtained the $L^p$-theory
of SPDEs on the whole space. Later, Kim \cite{K2004,K2008} established the $L^p$-theory of SPDEs on
the bounded domain. Using the Moser's iteration scheme, Denis et al. \cite{DMS} also obtained
the $L^p$-theory of SPDEs on the bounded domain. The other part is the Schauder estimates. Debussche et al.
\cite{DMH} proved that the solution of SPDEs is H\"{o}lder continuous in both time and space variables.
Du-Liu \cite{DL} established the $C^{2+\alpha}$-theory for SPDEs on the whole space. Hsu-Wang \cite{HW}
used stochastic De Giorgi iteration technique to prove that the solutions of SPDEs are almost surely
H\"{o}lder continuous in both space and time variables.

The above mentioned results about the regularity of the solutions of SPDEs belongs to the space
$L^p(\Omega;C^{\alpha,\beta}([0,T]\times G))$, where $G$ is a bounded domain in $\mathbb{R}^d$.
Now, there is an natural question, that is, can one get the H\"{o}lder estimate for the $p$-moment?
In other words, can we derive the estimate in $C^{\alpha,\beta}([0,T]\times G;L^p(\Omega))$?
We note that Du-Liu \cite{DL} obtained the $C^{2+\alpha}$-theory for SPDEs in
$C^{\alpha,\beta}([0,T]\times G);L^p(\Omega))$, where the Dini continuous is needed for the
stochastic term.
The method used in \cite{DL} is the Sobolev embedding theorem and the iteration technique under the condition that
the noise term satisfies Dini continuity. In this paper, we would like to consider a simple case, that is, the equation with
additive noise. We first derive the Morrey-Campanato estimates for the stochastic convolution operators and then, by
utilising the embedding theorem between Campanato space and H\"{o}lder space, we establish the norm of $C^{\theta,\theta/2}$.
As an application, we show that the solutions of partial differential equations driven by
Brownian motion or by L\'{e}vy noise are H\"{o}lder continuous in the both time and space variables on the whole space.
The novelty in our present paper is the approach we used is different from those in \cite{DMS,DL,HW}.
We would like to point out that by using the Morrey-Campanato estimates and the embedding theorem, the H\"{o}lder
estimates can be easily derived, and on the other hand our Morrey-Campanato estimates can be obtained by directly calculation,
thus our method is indeed simpler than others, see \cite[Lemma 4.3]{mLb} for the deterministic parabolic equations. Besides,
we establish the Schauder estimates for the solutions of partial differential equations driven by L\'{e}vy noise.

The paper is organized as follows. In Section 2, the next section, we set up our main results and present corresponding proofs.
Section 3, the final section, gives an application of our results.

Before ending up this section, let us introduce some notations. As usual, $\mathbb{R}^d$ stands for the $d$-dimensional Euclidean
space of points $x=(x_1,\cdots,x_d)$ with $|\cdot|$ being its usual Euclidean norm, and
$B_r(x):=\{y\in\mathbb{R}^d:|x-y|<r\}$ as well as $B_r:=B_r(0)$. We use $\mathbb{R}_+$ to denote the
set $\{x\in\mathbb{R}: x>0\}, a\wedge b:=\min\{a,b\}, a\vee b:=\max\{a,b\}$ and
$L^p:=L^p(\mathbb{R}^d)$. Finally, we write $N=N(a,b,\cdots)$ for a constant $N$ which depends on $a,b,\cdots$.

\section{ Main Results}
\setcounter{equation}{0}

We first recall some definitions and known results. Set, for $X=(t,x)\in\mathbb{R}\times\mathbb{R}^{d}$ and $Y=(s,y)\in\mathbb{R}\times\mathbb{R}^{d}$,
the following
   \bess
\delta(X,Y):=\max\left\{|x-y|,\,|t-s|^{\frac{1}{2}}\right\}.
   \eess
Let $Q_c(X)$ be the ball centered in $X=(t,x)$ and of radius $c$, i.e.,
   \bess
Q_c(X)&:=&\{Y=(s,y)\in\mathbb{R}\times\mathbb{R}^{d}:\,\delta(X,Y)<R\}\\
&=&(t-c^2,t+c^2)\times B_c(x).
  \eess
Fix $T\in(0,\infty]$ arbitrarily. Denote
   \bess
\mathcal {O}_T=(0,T)\times\mathbb{R}^d.
  \eess
Let $D$ be a bounded domain in $\mathbb{R}^{d+1}$ and for $X\in D$, $D(X,r):=D\cap Q_r(X)$ and
$d(D):=diam D$. We first introduce the definition of Campanato space.

\begin{defi}\lbl{d2.1} (Campanato Space) Let $p\geq1$ and $\theta\geq0$. $u$ belongs
to Campanato space $\mathscr{L}^{p,\theta}(D;\delta)$ if $u\in L^p(D)$ and
   \bess
[u]_{\mathscr{L}^{p,\theta}(D;\delta)}:=\left(\sup_{X\in D,d(D)\geq\rho>0}\frac{1}{|D(X,\rho)|^\theta}
\int_{D(X,\rho)}|u(Y)-u_{X,\rho}|^pdY\right)^{1/p}<\infty,
    \eess
and
    \bess
\|u\|_{\mathscr{L}^{p,\theta}(D;\delta)}:=\left(\|u\|_{L^p(D)}^p+[u]_{\mathscr{L}^{p,\theta}(D;\delta)}^p\right)^{1/p},
  \eess
where $|D(X,\rho)|$ stands for the Lebesgue measure of $D(X,\rho)$ and
   \bess
u_{X,\rho}=\frac{1}{|D(X,\rho)|}\int_{D(X,\rho)}u(Y)dY.
   \eess
\end{defi}

It is easy to verify that Campanato space is a Banach space and has the following property (see Appendix):
if $1\leq p\leq q<\infty$, $(\theta-p)/p\leq(\sigma-p)/q$, it holds that
   \bess
\mathscr{L}^{q,\sigma}(D;\delta)\subset\mathscr{L}^{p,\theta}(D;\delta).
  \eess
Next we recall the definition of H\"{o}lder space.

\begin{defi}\lbl{d2.2} (H\"{o}lder Space) Let $0<\alpha\leq1$. $u$ belongs
to H\"{o}lder space $C^{\alpha}(\bar D;\delta)$ if $u$ satisfies
   \bess
[u]_{C^{\alpha}(\bar D;\delta)}:=\sup_{X\in D,d(D)\geq\rho>0}\frac{|u(X)-u(Y)|}{\delta(X,Y)^\alpha}
<\infty,
    \eess
and
    \bess
\|u\|_{C^{\alpha}(\bar D;\delta)}:=\sup_D|u|+[u]_{C^{\alpha}(\bar D;\delta)}.
  \eess
\end{defi}

\begin{defi}\lbl{d2.3} Let $D\subset\mathbb{R}^{d+1}$. Domain $D$ is called $A$-type if there exists a constant $A>0$ such that
$\forall X\in D$, $0<\rho\leq diam D$, it holds that
   \bess
|D(X,\rho)|\geq A|Q_\rho(X)|.
  \eess
\end{defi}

Comparing with the two space, we have the following relations.
\begin{prop}\lbl{p2.1} Assume that $D$ is an $A$-type bounded domain.
Then we have the following relation: when $1<\theta\leq1+\frac{p}{d+2}$ and $p\geq1$,
   \bess
\mathscr{L}^{p,\theta}(D;\delta)\cong C^{\alpha}(\bar D;\delta),
  \eess
where
   \bess
\alpha=\frac{(d+2)(\theta-1)}{p},
   \eess
where $d$ is the dimension of the space.
\end{prop}

Here $A\cong B$ means that both $A\subseteq B$ and $B\subseteq A$ hold. We will
obtain the Campanato estimates under some assumptions on the kernel $K$.
Noting that
   \bess
&&\left(\sup_{X\in D,d(D)\geq\rho>0}\frac{1}{|D(X,\rho)|^\theta}
\int_{D(X,\rho)}|u(Y)-u_{X,\rho}|^pdY\right)^{1/p}\\
&=&\left(\sup_{X\in D,d(D)\geq\rho>0}\frac{1}{|D(X,\rho)|^\theta}
\int_{D(X,\rho)}\Big|u(Y)-\frac{1}{|D(X,\rho)|}\int_{D(X,\rho)}u(Z)dZ\Big|^pdY\right)^{1/p}\\
&\leq&\left(\sup_{X\in D,d(D)\geq\rho>0}\frac{1}{|D(X,\rho)|^{1+\theta}}
\int_{D(X,\rho)}\int_{D(X,\rho)}|u(Y)-u(Z)|^pdZdY\right)^{1/p},
  \eess
so the definition of semi-norm of the Campanato space can be replaced by the
above inequality. We also remark that in order to get the H\"{o}lder estimate, the
range of $\theta$ must be lager than $1$.

Now, we talk about two spaces $L^p(\Omega;\mathscr{L}^{p,\theta}(D;\delta))$ and
$\mathscr{L}^{p,\theta}((D;\delta);L^p(\Omega))$. If we want to prove
$u\in L^p(\Omega;\mathscr{L}^{p,\theta}(D;\delta))$, that is,
   \bess
\mathbb{E}[u]^p_{\mathscr{L}^{p,\theta}(D;\delta)}:=\mathbb{E}\sup_{D}\frac{1}{|D|^{1+\theta}}
\int_D\int_D|h(t,x)-h(s,y)|^pdtdxdsdy<\infty,
    \eess
our first idea is to prove the two maps $\mathbb{E}$ and $\sup_{t,x}$ can be interchanged.
Unfortunately, it is hard to give a sufficient condition to assure the above idea holds.
The second idea is to prove the norm of $u$ in $\mathscr{L}^{p,\theta}(D;\delta)$ is bounded
almost surely. The two ideas is hard to come true.
And thus we must adjust our idea.
We also remark that the mean of the space
$\mathscr{L}^{p,\theta}((D;\delta);L^p(\Omega))$ is that we call
   \bess
u\in\mathscr{L}^{p,\theta}((D;\delta);L^p(\Omega)),\ \ \ {\rm if}\ \ \
\|u\|_{L^p(\Omega)}\in \mathscr{L}^{p,\theta}(D;\delta).
   \eess
In other words, the following norm is finite
   \bess
&&[\|u\|_{L^p(\Omega)}]^p_{\mathscr{L}^{p,\theta}(D;\delta)}\\
&:=&\sup_{D}\frac{1}{|D|^{1+\theta}}
\int_D\int_D\Big|\|u\|_{L^p(\Omega)}(t,x)-\|u\|_{L^p(\Omega)}(s,y)\Big|^pdtdxdsdy\\
&<&\infty.
    \eess
Using triangular inequality, we have
    \bess
&&[\|u\|_{L^p(\Omega)}]^p_{\mathscr{L}^{p,\theta}(D;\delta)}\\
&\leq&\sup_{D}\frac{1}{|D|^{1+\theta}}
\int_D\int_D\|u(t,x)-u(s,y)\|^p_{L^p(\Omega)}dtdxdsdy\\
&=&\sup_{D}\frac{1}{|D|^{1+\theta}}\mathbb{E}
\int_D\int_D|u(t,x)-u(s,y)|^pdtdxdsdy.
    \eess
Thus we only need to prove that
   \bess
\sup_{D}\frac{1}{|D|^{1+\theta}}\mathbb{E}
\int_D\int_D|u(t,x)-u(s,y)|^pdtdxdsdy<\infty.
   \eess

\subsection{Brownian Motion Case}

Let $(\Omega,\mathcal {F},\mathbb{P})$ be a complete probability space endowed with
$\{\mathcal {F}_t\}_{t\in[0,T]}$, a filtration on $\Omega$ containing
all $P$-null subsets of $\Omega$. Let $W_t$ be a one-dimensional $\{\mathcal {F}_t\}_{t\in[0,T]}$-adapted Wiener processes
defined on $(\Omega,\mathcal {F},\mathbb{P})$.

For a measurable function $h$ on $\Omega\times\mathcal {O}_T$, we define the
stochastic Campanato quasi-norm of $h$ on $\Omega\times\mathcal {O}_T$
as follows:
   \bess
[h]^p_{\mathscr{L}^{p,\theta}((Q;\delta);L^p(\Omega))}:=\sup_{Q}\frac{1}{|Q|^{1+\theta}}\mathbb{E}
\int_Q\int_Q|h(t,x)-h(s,y)|^pdtdxdsdy,
    \eess
where the sup is taken over all $Q=D\cap Q_c$ of the type
   \bess
Q_c(t_0,x_0):=(t_0-c^2,t_0+c^2)\times B_c(x_0)\subset\mathcal {O}_T,\ \ c>0,t_0>0.
  \eess
It is remarked that when $\theta=1$, this is equivalent to the classical BMO semi-norm which is
introduced in John-Nirenberg \cite{JN}. If the stochastic Campanato quasi-norm of $h$ is
finite, we then say that $h$ belongs to the space $\mathscr{L}^{p,\theta}((Q;\delta);L^p(\Omega))$.

We first consider the Brownian motion case. Given a deterministic kernel
$K:\mathbb{R}\times\mathbb{R}^d\to\mathbb{R}$, we denote for any no-random (i.e., not randomly dependent) $g:\mathbb{R}\times\mathbb{R}^d\to\mathbb{R}$
the following stochastic convolution
  \bes
\mathcal {K}g(t,x):=\int_0^t\int_{\mathbb{R}^d}K(t-r,y)g(r,x-y)dydW(r).
  \lbl{2.1}\ees

Then we have the following result.
  \begin{theo}\lbl{t2.1}
Let $D$ be an $A$-type bounded domain in $\mathbb{R}^{d+1}$ such that $\bar D\subset \mathcal {O}_T$.
Suppose that $g\in C^\beta(\mathbb{R}_+\times\mathbb{R}^d)$, $0<\beta<1$, is a non-random function
  and $g(0,0)=0$.
Assume that there exists positive constants $\gamma_i$ ($i=1,2$) such that
the non-random kernel function satisfies that for any $t\in(0,T]$
  \bes
&&\int_0^s\left(\int_{\mathbb{R}^d}|K(t-r,z)-K(s-r,z)|(1+|z|^\beta)dz\right)^2dr\leq N(T,\beta)(t-s)^{\gamma_1},\lbl{2.2}\\
&&\int_0^s\left(\int_{\mathbb{R}^d}|K(s-r,z)|dz\right)^2dr\leq N_0,\lbl{2.3}\\
&& \int_s^t\left(\int_{\mathbb{R}^d}|K(t-r,z)|(1+|z|^\beta) dz\right)^2dr\leq N(T,\beta)(t-s)^{\gamma_2},
   \lbl{2.4} \ees
where $N_0$ is a positive constant. Then we have, for
$p\geq1$ and $\beta<\gamma$,
    \bess
[\mathcal {K}g]_{\mathscr{L}^{p,\theta}((D;\delta);L^p(\Omega))}\leq N(N_0,\beta,T,d,p),
    \eess
where $\theta=1+\frac{\gamma p}{d+2}$ and $\gamma=\min\{\gamma_1,\gamma_2,\beta\}$.
\end{theo}

{\bf Proof.} Let $(t_0,x_0)\in D\subset\mathcal {O}_T$ and
  \bess
Q_c(t_0,x_0)=(t_0-c^2,t_0+c^2)\times B_c(x_0).
   \eess
Then set $C_1:=diam D$, we have $\bar D\subset Q_{C_1}(t_0,x_0)$.
Denote $Q:=D\cap Q_c(t_0,x_0)$.

Set $t>s$. By the BDG inequality, we have
  \bess
&&\mathbb{E}\int_Q\int_Q|\mathcal {K}g(t,x)-\mathcal {K}g(s,y)|^pdtdxdsdy\\
&=&\mathbb{E}\int_Q\int_Q\Big|\int_0^t\int_{\mathbb{R}^d}K(t-r,z)g(r,x-z)dzdW(r)\\
&&
-\int_0^s\int_{\mathbb{R}^d}K(s-r,z)g(r,y-z)dzdW(r)\Big|^pdtdxdsdy\\
&\leq&2^{p-1}\mathbb{E}\int_Q\int_Q\Big|\int_0^s\int_{\mathbb{R}^d}(K(t-r,z)-K(s-r,z))
g(r,x-z)dzdW(r)\Big|^p \\
&&+2^{p-1}\mathbb{E}\int_Q\int_Q\Big|\int_0^s\int_{\mathbb{R}^d}K(s-r,z)
(g(r,x-z)-g(r,y-z))dzdW(r)\Big|^p \\
&&+2^{p-1}\mathbb{E}\int_Q\int_Q\Big|\int_s^t\int_{\mathbb{R}^d}K(t-r,z)
g(r,x-z)dzdW(r)\Big|^pdtdxdsdy\\
&\leq&N(p)\int_Q\int_Q
\left(\int_0^s|\int_{\mathbb{R}^d}|K(t-r,z)-K(s-r,z)||g(r,x-z)|dz|^2
dr\right)^{\frac{p}{2}} \\
&&
+N(p)\int_Q\int_Q
\left(\int_0^s|\int_{\mathbb{R}^d}|K(s-r,z)||g(r,x-z)-g(r,y-z)|dz|^2
dr\right)^{\frac{p}{2}} \\
&&
+N(p)\int_Q\int_Q\left(\int_s^t|
\int_{\mathbb{R}^d}K(t-r,z)g(r,x-z)dz|^2dr\right)^{\frac{p}{2}}\\
&=:&I_1+I_2+I_3.
  \eess

Estimate of $I_1$. By using the H\"{o}lder continuous of $g$, i.e.,
   \bess
|g(r,x-z)-g(0,0)|&\leq& C_g\max\left\{r^{\frac{1}{2}},|x-z|\right\}^\beta\\
&\leq&N(g,\beta)(T^{\frac{\beta}{2}}+|x-x_0|^\beta+|x_0|^\beta+|z|^\beta)\\
&\leq&N(g,\beta)(T^{\frac{\beta}{2}}+c^\beta+|x_0|^\beta+|z|^\beta),
   \eess
and (\ref{2.2}), we have
   \bess
I_1&=&N(p)\int_Q\int_Q
\left(\int_0^s|\int_{\mathbb{R}^d}|K(t-r,z)-K(s-r,z)||g(r,x-z)|dz|^2
dr\right)^{\frac{p}{2}}\\
&\leq&N(p,\beta)\int_Q\int_Q
\left(\int_0^s|\int_{\mathbb{R}^d}|K(t-r,z)-K(s-r,z)|(T^{\frac{\beta}{2}}+c^\beta+|x_0|^\beta+|z|^\beta)dz|^2
dr\right)^{\frac{p}{2}}\\
&\leq&N(p,\beta,T,x_0)\int_Q\int_Q
\left(\int_0^s|\int_{\mathbb{R}^d}|K(t-r,z)-K(s-r,z)|(1+|z|^\beta)dz|^2
dr\right)^{\frac{p}{2}}\\
&&+c^\beta N(p,\beta)\int_Q\int_Q
\left(\int_0^s|\int_{\mathbb{R}^d}|K(t-r,z)-K(s-r,z)
dr\right)^{\frac{p}{2}}\\
&\leq&N(p,\beta,T,x_0)(1+|x-y|^{\beta p})(t-s)^{{\frac{\gamma_1p}{2}}}|Q|^2.
   \eess

The condition (\ref{2.3}) and
   \bess
|g(r,x-z)-g(r,y-z)|\leq C_g|x-y|^\beta
   \eess
imply the following derivation
   \bess
I_2&=&N(p)\int_Q\int_Q\left(\int_0^s|\int_{\mathbb{R}^d}|K(s-r,z)||g(r,x-z)-g(r,y-z)|dz|^2dr\right)^{\frac{p}{2}}\\
&\leq&N(p,g)\int_Q\int_Q\left(\int_0^s|\int_{\mathbb{R}^d}|K(r,z)||x-y|^\beta dz|^2dr\right)^{\frac{p}{2}}\\
&\leq& N(N_0,p,g,\beta)|x-y|^{\beta p}|Q|^2.
   \eess

Estimate of $I_3$. By using the property $g(0,0)=0$ and (\ref{2.4}), we get
   \bess
I_3&=&N(p)\int_Q\int_Q\left(\int_s^t|
\int_{\mathbb{R}^d}K(t-r,z)g(r,x-z)dz|^2dr\right)^{\frac{p}{2}}\\
&\leq&\int_Q\int_Q\left(\int_s^t\Big|
\int_{\mathbb{R}^d}|K(r,z)|
(T+|x-x_0|^\beta+|x_0|^\beta+|z|^\beta) dz\Big|^2dr\right)^{\frac{p}{2}}\\
&\leq&N(p,T,x_0,\beta)\int_Q\int_Q
\left(\int_s^t\Big|\int_{\mathbb{R}^d}|K(t-r,z)|(1+|z|^\beta) dz
\Big|^2dr\right)^{\frac{p}{2}}\\
&&+N(p,T,\beta)|x-y|^{\beta p}\int_Q\int_Q
\left(\int_s^t\Big|\int_{\mathbb{R}^d}|K(t-r,z)| dz
\Big|^2dr\right)^{\frac{p}{2}}\\
&\leq&N(p,T,x_0,\beta)|Q|^2(t-s)^{\frac{\gamma_2p}{2}}(1+|x-y|^{\beta p}).
   \eess

Noting that $(t,x)\in Q_c$ and $(s,y)\in Q_c$, we have
   \bess
0\leq t-s\leq 2c^2\ \ \ {\rm and}\ \ \ |x-y|\leq|x-x_0|+|y-x_0|\leq2c.
   \eess
Using the above inequality and the properties of $A$-type domain, we deduce
  \bess
  I_1&\leq&N(p,T,\beta,x_0)(1+c^{\beta p})c^{\gamma_1p}|Q|^2;\\
I_2&\leq& N(N_0,p,g,\beta)c^{\beta p}|Q|^2;\\
I_3&\leq&N(p,T,x_0,\beta)|Q|^2c^{\gamma_2p}(1+c^{\beta p}).
   \eess
Combining the estimates of $I_1,I_2$ and $I_3$, we get
 \bess
&&\mathbb{E}\int_Q\int_Q|u(t,x)-u(s,y)|^pdtdxdsdy\\
&\leq&N(\beta,N_0,T,p)|Q|^2(c^{\beta p}+1)(c^{\beta p}+c^{\gamma_1p}+c^{\gamma_2p}).
   \eess
Since $D$ is a $A$-type bounded domain, we have $c\leq diam D$ and
   \bess
A|Q_c(t_0,x_0)|\leq |Q|\leq |Q_c(t_0,x_0)|.
   \eess
We remark that $|Q_c(t_0,x_0)|=Nc^{d+2}$ and $0<\beta\leq1$, where $N$ is a positive constant which does not depend on
$c$. Noting that $Q\subset Q_{C_1}$, we have
  \bess
&&\mathbb{E}\int_Q\int_Q|u(t,x)-u(s,y)|^pdtdxdsdy\\
&\leq&N(\beta,N_0,C_1,d,T)|Q|^{2+\frac{\gamma p}{d+2}},
   \eess
where $\gamma=\min\{\gamma_1,\gamma_2,\beta\}$, which yield that
      \bess
&&[\mathcal {K}g]_{\mathscr{L}^{p,\theta}(D;\delta)}\\
&=&\sup_{Q}\frac{1}{|Q|^{1+\theta}}\mathbb{E}\int_Q\int_Q|\mathcal {K}g(t,x)-\mathcal {K}g(s,y)|^pdtdxdsdy\\
&\leq& N(\beta,N_0,T,d,p),
    \eess
where $\theta=1+\frac{\gamma p}{d+2}$. The proof of Theorem \ref{t2.1} is complete. $\Box$

Theorem \ref{t2.1} shows that $\mathcal {K}g(t,x)\in \mathscr{L}^{p,\theta}((Q;\delta);L^p(\Omega))$.
That is, $\|\mathcal {K}g\|_{L^p(\Omega)}\in \mathscr{L}^{p,\theta}(Q;\delta)$.
Applying the result of Proposition \ref{p2.1}, we have the following result.
\begin{col}\lbl{c2.1}
Assume all the assumptions in Theorem \ref{t2.1} hold, then
   \bess
\mathcal {K}g(t,x)\in C^{\gamma}((\bar D;\delta);L^p(\Omega)).
  \eess
\end{col}

\begin{remark}\lbl{r2.1} 1. It follows from Theorem \ref{t2.1} and Corollary \ref{c2.1} that
$\mathcal {K}g(t,x)\in C^{\gamma}((\bar D;\delta);L^p(\Omega))$ and $\gamma=\min\{\gamma_1,\gamma_2,\beta\}$ if $g\in C^\beta(\mathbb{R}_+\times\mathbb{R}^d)$
and $g(0,0)=0$. For special kernel, we can let $\gamma=\beta$, see Theorem \ref{t3.1}. That is to say, the regularity of $\mathcal {K}g(t,x)$ depends
heavily on the noise term $g$.

2. It is easy to prove that if $g\in C^{k+\beta,\beta/2}(\mathbb{R}_+\times\mathbb{R}^d)$
and $\nabla^kg(0,0)=0$, then $\mathcal {K}g(t,x)\in C^{k+\gamma,\gamma/2}(\bar D;\delta)$ under
the assumptions of Theorem \ref{t2.1}.
Here $g\in C^{k+\beta,\beta/2}(\mathbb{R}_+\times\mathbb{R}^d)$ denotes that the $k$-order of $g$ w.r.t
spatial variable belongs to $C^\beta$, and that $g$ w.r.t time variable belongs to $C^{\beta/2}$.

3. The regularity w.r.t time variable can not be improved because of the fact that the regularity of Brownian motion w.r.t time variable is $C^{\frac{1}{2}-}$.

4. If the kernel function $K$ is random, the similar result also holds. The constant $N$ in
Theorem \ref{t2.1} depending on the choice of $x_0$ can be removed provided that
   \bess
\mathbb{E}\left[\|g\|^{p_0}_{L^\infty(\mathcal {O}_T)}\right]<\infty,
  \eess
where $p_0\geq1$ and $1\leq p\leq p_0$.

5. The method used in Theorem \ref{t2.1} is similar to that in \cite{mLb} for the
interior Schauder estimate, see \cite[Lemma 4.3]{mLb}.
\end{remark}

In Theorem \ref{t2.1}, the noise term $g$ depends on the times and spatial variables.
A natural question is: if $g$ does not depend on the time $t$, the
result of Theorem \ref{t2.1} will also hold ? Next, we answer this question.

  \begin{theo}\lbl{t2.2} Suppose that $g\in C^\beta(\mathbb{R}^{d+1})$, $0<\beta<1$
  and $g(0)=0$. Assume further that (\ref{2.2})-(\ref{2.4}) hold.
Let $D$ be a $A$-type bounded domain in $\mathbb{R}^{d+1}$ such that $\bar D\subset \mathcal {O}_T$ Then we have, for
$p\geq1$,
    \bess
[\mathcal {K}g]_{\mathscr{L}^{p,\theta}(D;\delta)}\leq N(N_0,\beta,T,d,p),
    \eess
where $\theta=1+\frac{\gamma p}{d+2}$ and $\gamma=\min\{\gamma_1,\gamma_2,\beta\}$.
\end{theo}

{\bf Proof.} The definition of $Q$ is the same as in the proof of Theorem \ref{t2.1}. Fix $t>s$.
The BDG inequality implies that
  \bess
&&\mathbb{E}\int_Q\int_Q|\mathcal {K}g(t,x)-\mathcal {K}g(s,y)|^pdtdxdsdy\\
&=&\mathbb{E}\int_Q\int_Q\Big|\int_0^t\int_{\mathbb{R}^d}K(t-r,z)g(x-z)dzdW(r)\\
&&
-\int_0^s\int_{\mathbb{R}^d}K(s-r,z)g(y-z)dzdW(r)\Big|^pdtdxdsdy\\
&\leq&N(p)\mathbb{E}\int_Q\int_Q\Big|\int_0^s\int_{\mathbb{R}^d}(K(t-r,z)-K(s-r,z))g(y-z)dzdW(r)\Big|^p \\
&&+N(p)\mathbb{E}\int_Q\int_Q\Big|\int_0^s\int_{\mathbb{R}^d}K(t-r,z)(g(x-z)-g(y-z))dzdW(r)\Big|^p \\
&&
+N(p)\mathbb{E}\int_Q\int_Q\Big|\int_s^t\int_{\mathbb{R}^d}K(t-r,z)g(x-z)dzdW(r)\Big|^pdtdxdsdy\\
&\leq&N(p)\mathbb{E}\int_Q\int_Q\left(\int_0^s|\int_{\mathbb{R}^d}(K(t-r,z)-K(s-r,z))g(y-z)dz|^2dr\right)^{\frac{p}{2}} \\
&&+N(p)\mathbb{E}\int_Q\int_Q\left(\int_0^s|\int_{\mathbb{R}^d}K(t-r,z)(g(x-z)-g(y-z))dz|^2dr\right)^{\frac{p}{2}} \\
&&+N(p)
\mathbb{E}\int_Q\int_Q\left(\int_s^t|\int_{\mathbb{R}^d}K(t-r,z)g(x-z)dz|^2dr\right)^{\frac{p}{2}}\\
&=:&I_1+I_2+I_3.
  \eess
Noting again that $(t,x)\in Q_c$ and $(s,y)\in Q_c$, we have
   \bess
0\leq t-s\leq 2c^2\ \ \ {\rm and}\ \ \ |x-y|\leq|x-x_0|+|y-x_0|\leq2c.
   \eess
The H\"{o}lder continuous of $g$ and (\ref{2.2})-(\ref{2.4}) give that
   \bess
I_1+I_2+I_3&\leq& N(N_0,p,\beta,T,d)|Q|^2(c^{\beta p}+c^{\gamma_1p}+c^{\gamma_2p}),
   \eess
which implies the desired result. The proof is complete. $\Box$
\begin{remark}\lbl{r2.2} By using Proposition \ref{p2.1}, one can get
$\mathcal {K}g(t,x)\in C^{\gamma}((\bar D;\delta);L^p(\Omega))$. In particular, taking
$g=constant$, we have the regularity of time variable is $C^{\frac{1}{2}-}$ and
the regularity of spatial variable is $C^\infty$.
  \end{remark}

\subsection{L\'{e}vy Noise Case}
Let $(\Omega,\mathcal {F},\mathbb{F},\mathbb{P})$ be a complete probability space such that
$\{\mathcal {F}_t\}_{t\in[0,T]}$ is a filtration on $\Omega$ containing
all $P$-null subsets of $\Omega$ and $\mathbb{F}$ be the predictable $\sigma$-algebra associated with the filtration $\{\mathcal {F}_t\}_{t\in[0,T]}$. We are given a $\sigma$-finite measure space $(Z,\mathcal {Z},\nu)$ and a Poisson
random measure $\mu$ on $[0,T]\times Z$, defined on the stochastic basis. The compensator of
$\mu$ is Leb$\otimes\nu$, and the compensated martingale measure $\tilde N:=\mu-Leb\otimes\nu$.

In this subsection, we consider the stochastic singular integral operator
    \bes
\mathcal {G}g(t,x)&=&\int_0^t\int_ZK(t,s,\cdot)\ast g(s,\cdot,z)(x)\tilde N(dz,ds)\nm\\
&=&\int_0^t\int_Z\int_{\mathbb{R}^d}K(t-s,x-y)g(s,y,z)dy\tilde N(dz,ds)
   \lbl{2.5}\ees
for $\mathbb{F}$-predictable processes $g:[0,T]\times\mathbb{R}^d\times Z\times\Omega\to\mathbb{R}$.
For simplicity, we assume that the kernel function is deterministic.
We first recall the Kunita's first inequality.

\begin{defi}\lbl{d2.4}(Kunita's first inequality \cite[Theorem 4.4.23]{Ap}) For any $p\geq2$,
there exists $N(p)>0$ such that
   \bes
\mathbb{E}\left(\sup_{0\leq t\leq T}|I(t)|^p\right)&\leq&N(p)\left\{
\mathbb{E}\left[\left(\int_0^T\int_Z|H(t,z)|^2\nu(dz)dt\right)^{p/2}\right]\right.\nm\\
&&\left.+\mathbb{E}\left[\int_0^T\int_Z|H(t,z)|^p\nu(dz)dt\right]\right\},
   \lbl{2.6}\ees
where $H\in \mathcal {P}_2(t,E)$ and
    \bess
I(t)=\int_0^t\int_ZH(s,z)\tilde N(dz,ds).
   \eess
$\mathcal {P}_2(T,E)$ denotes the set of all equivalence classes of mappings
$F:[0,T]\times E\times\Omega\to\mathbb{R}$ which coincide almost everywhere with respect to
$\rho\times P$ and which satisfy the following conditions (see Page 225 of \cite{Ap})

(i) $F$ is $\mathbb{F}$-predictable;

(ii) $P\left(\int_0^T\int_E|F(t,x)|^2\rho(dt,dx)<\infty\right)=1$.
   \end{defi}

Now we are in the position to show our main result.
  \begin{theo}\lbl{t2.3} Let $g_1:Z\times\Omega\to\mathbb{R}$ be measurable and fulfil the following
         \bess
\mathbb{E}\left[\left(\int_Z|g_1(z)|^2\nu(dz)\right)^{p_0/2}+\int_Z|g_1(z)|^{p_0}\nu(dz)\right]<\infty
  \eess
for some constant $p_0>2$. Suppose that the function $g$ satisfies that
    \bes
  |g(t,x,z)-g(s,y,z)|\leq C_g\max\left\{(t-s)^{\frac{1}{2}},|x-y|\right\}^\beta g_1(z),\ \ \
  {\rm for \ all}\ \ z\in Z,\ \ a.s.,
     \lbl{2.7}\ees
and $g(0,0,z)=0$ uniformly for $z\in Z$ almost surely.  Assume further that there exist positive constants $\gamma_i$ ($i=1,2$) such that
the non-random kernel function satisfies that for any $t\in(0,T]$,
   \bess
&&\int_0^s\left(\int_{\mathbb{R}^d}|K(t-r,z)-K(s-r,z)|(1+|z|^\beta)dz\right)^pdr\leq N(T,\beta)(t-s)^{\frac{\gamma_1p}{2}}, \\
&&\int_0^s\left(\int_{\mathbb{R}^d}|K(s-r,z)|dz\right)^pdr\leq N_0, \\
&& \int_s^t\left(\int_{\mathbb{R}^d}|K(t-r,z)|(1+|z|^\beta) dz\right)^pdr\leq N(T,\beta)(t-s)^{\frac{\gamma_2p}{2}},
     \eess
where $N_0$ is a positive constant.
Let $D$ be an $A$-type bounded domain in $\mathbb{R}^{d+1}$ such that $\bar D\subset \mathcal {O}_T$. Then we have, for
$2\leq p\leq p_0$ and $\beta<\alpha$,
    \bess
[\mathcal {K}g(t,x)]_{\mathscr{L}^{p,\theta}(D;\delta)}\leq N(N_0,\beta,T,d,p),
    \eess
where $\theta=1+\frac{\gamma p}{d+2}$ and $\gamma=\min\{\gamma_1,\gamma_2,\beta\}$.
\end{theo}

{\bf Proof.} Similar to the proof of Theorem \ref{t2.1} and using the inequality (\ref{2.6}) we first have the following estimates.
     \bess
&&\mathbb{E}|\mathcal {G}g(t,x)-\mathcal {G}g(s,y)|^p\\
&=&\mathbb{E}\left[\Big|\int_0^t\int_Z\int_{\mathbb{R}^d}K(t-r,\xi)g(r,x-\xi,z)d\xi\tilde N(dz,dr)\right.\\
&&\qquad\left.
-\int_0^s\int_Z\int_{\mathbb{R}^d}K(s-r,\xi)g(r,y-\xi,z)d\xi\tilde N(dz,dr)\Big|^p\right]\\
&\leq&N(p)\mathbb{E}\left[\Big|\int_0^s\int_Z\int_{\mathbb{R}^d}[K(t-r,\xi)-K(s-r,\xi)]g(r,x-\xi,z)d\xi\tilde N(dz,dr)\right.\\
&&\qquad+\int_0^s\int_Z\int_{\mathbb{R}^d}K(s-r,\xi)[g(r,x-\xi,z)-g(r,y-\xi,z)]d\xi\tilde N(dz,dr)\\
&&\qquad\left.
+\int_s^t\int_Z\int_{\mathbb{R}^d}K(t-r,\xi)g(r,x-\xi,z)d\xi\tilde N(dz,dr)\Big|^p\right]\\
&\leq&N(p)\mathbb{E}\left[\left(\int_s^t\int_Z|\int_{\mathbb{R}^d}K(t-r,\xi)g(r,x-\xi,z)d\xi|^2\nu(dz)dr\right)^{p/2}\right]\nm\\
&& +N(p)\mathbb{E}\left[\int_s^t\int_Z|\int_{\mathbb{R}^d}K(t-r,\xi)g(r,x-\xi,z)d\xi|^p\nu(dz)dr\right] \\
&&+N(p)\mathbb{E}\left[\left(\int_0^s\int_Z|\int_{\mathbb{R}^d}K(s-r,\xi)[g(r,x-\xi,z)-g(r,y-\xi,z)]d\xi|^2\nu(dz)dr\right)^{p/2}\right]\nm\\
&& +N(p)\mathbb{E}\left[\int_0^s\int_Z|\int_{\mathbb{R}^d}K(s-r,\xi)[g(r,x-\xi,z)-g(r,y-\xi,z)]d\xi|^p\nu(dz)dr\right]\\
&&+N(p)\mathbb{E}\left[\left(\int_0^s\int_Z|\int_{\mathbb{R}^d}[K(t-r,\xi)-K(s-r,\xi)]g(r,x-\xi,z)d\xi|^2\nu(dz)dr\right)^{p/2}\right]\nm\\
&& +N(p)\mathbb{E}\left[\int_0^s\int_Z|\int_{\mathbb{R}^d}[K(t-r,\xi)-K(s-r,\xi)]g(r,x-\xi,z)d\xi|^p\nu(dz)dr\right].
  \eess
By using (\ref{2.7}) and $g(0,0,z)=0$ uniformly for $z\in Z$ almost surely, we have that the above inequality
is smaller than or equal to
   \bes
&&N\mathbb{E}\left[\left(\int_s^t\int_Zg_1(z)^2\Big|\int_{\mathbb{R}^d}|K(t-r,\xi)|(|x-x_0|^\beta+|x_0-\xi|^\beta)
d\xi\Big|^2\nu(dz)dr\right)^{p/2}\right]\nm\\
&& +N\mathbb{E}\left[\int_s^t\int_Z|g_1(z)|^p\Big|\int_{\mathbb{R}^d}|K(t-r,\xi)|(|x-x_0|^\beta+|x_0-\xi|^\beta) d\xi\Big|^p\nu(dz)dr\right]\nm \\
&&+N\mathbb{E}\left[\left(\int_0^s\int_Zg_1(z)^2\Big|\int_{\mathbb{R}^d}|K(s-r,\xi)||x-y|^\beta
d\xi\Big|^2\nu(dz)dr\right)^{p/2}\right]\nm\\
&& +N\mathbb{E}\left[\int_0^s\int_Z|g_1(z)|^p\Big|\int_{\mathbb{R}^d}|K(s-r,\xi)||x-y|^\beta d\xi\Big|^p\nu(dz)dr\right]\nm\\
&&+N\mathbb{E}\left[\left(\int_0^s\int_Zg_1(z)^2|\int_{\mathbb{R}^d}|K(t-r,\xi)-K(s-r,\xi)|(|x-x_0|^\beta+|x_0-\xi|^\beta)d\xi|^2\nu(dz)dr\right)^{p/2}\right]\nm\\
&& +N\mathbb{E}\left[\int_0^s\int_Z|g_1(z)|^p|\int_{\mathbb{R}^d}|K(t-r,\xi)-K(s-r,\xi)|(|x-x_0|^\beta+|x_0-\xi|^\beta)d\xi|^p\nu(dz)dr\right].
  \lbl{2.8}\ees
Following the proof of Theorem \ref{t2.1}, we have
   \bess
0\leq t-s\leq 2c^2\ \ \ {\rm and}\ \ \ |x-y|\leq|x-x_0|+|y-x_0|\leq2c.
   \eess
Thus (\ref{2.8}) yields that
   \bess
&&\mathbb{E}|\mathcal {G}g(t,x)-\mathcal {G}g(s,y)|^p\\
&\leq&N(p,T,|x_0|)(1+c^{\beta p})\mathbb{E}\left[\left(\int_s^t\int_Zg_1(z)^2\Big|\int_{\mathbb{R}^d}|K(t-r,\xi)|(1+|\xi|^\beta)
d\xi\Big|^2\nu(dz)dr\right)^{p/2}\right]\nm\\
&&+N(p,T)(1+c^{\beta p})\mathbb{E}\left[\int_s^t\int_Z|g_1(z)|^p\Big|\int_{\mathbb{R}^d}|K(r,\xi)|(1+|\xi|^\beta)
d\xi\Big|^p\nu(dz)dr\right]\nm\\
&&+N(p,T)c^{\beta p}\mathbb{E}\left[\left(\int_0^s\int_Zg_1(z)^2\Big|\int_{\mathbb{R}^d}|K(r,\xi)|d\xi\Big|^2\nu(dz)dr\right)^{p/2}\right]\nm\\
&& +N(p,T)c^{\beta p}\mathbb{E}\left[\int_0^s\int_Z|g_1(z)|^p\Big|\int_{\mathbb{R}^d}|K(r,\xi)| d\xi\Big|^p\nu(dz)dr\right]\\
&&+N(p,T)(1+c^{\beta p})\mathbb{E}\left[\left(\int_0^s\int_Zg_1(z)^2\Big|\int_{\mathbb{R}^d}|K(t-r,\xi)-K(s-r,\xi)|(1+|\xi|^\beta)d\xi\Big|^2\nu(dz)dr\right)^{p/2}\right]\nm\\
&& +N(p,T)(1+c^{\beta p})\mathbb{E}\left[\int_0^s\int_Z|g_1(z)|^p\Big|\int_{\mathbb{R}^d}|K(t-r,\xi)-K(s-r,\xi)|(1+|\xi|^\beta) d\xi\Big|^p\nu(dz)dr\right]\\
&\leq& N(\beta,p,T,N_0)[1+c^{(1-\beta)p}](c^{\gamma_1 p}+c^{\gamma_2 p}+c^{\beta p}).
   \eess
Similar to the proof of Theorem \ref{t2.1}, by using the properties of $A$-type domain, one can complete the proof of
Theorem \ref{t2.3}.  $\Box$

\begin{col}\lbl{c2.2}
Assume all the assumptions in Theorem \ref{t2.3} hold, then
   \bess
\mathcal {G}g(t,x)\in C^{\gamma}((\bar D;\delta);L^p(\Omega)).
  \eess
\end{col}

\begin{remark}\lbl{l2.3} In Theorem \ref{t2.3}, both indices $\gamma_i, i=1,2$, depend on
the parameter $p$. On the other hand, we notice that when $p=2$, the two indices $\gamma_i, i=1,2$ will coincide with those in Theorem \ref{t2.1}. It then follows from Proposition \ref{p2.1} that $p\geq1$ is necessary and hence we can let $p=2$. Moreover, $\gamma$ will be largest in case $p=2$.
 \end{remark}

\section{Applications}
\setcounter{equation}{0}

In this section, applying Theorems \ref{t2.1}, \ref{t2.2} and \ref{t2.3},
we give some examples.

\subsection{Application to Parabolic Equations Driven by Brownian Motion}
In this subsection, We first consider the following stochastic
parabolic equations
    \bes
\left\{\begin{array}{lll}
du(t,x)=(\Delta u+divB(u)+c(t,x)u+f(t,x))dt+g(t,x)dW(t),\ \ t>0,\,x\in\mathbb{R}^d,\\
u(0,x)=u_0(x),\ \  \,x\in\mathbb{R}^d.
  \end{array}\right.
  \lbl{3.1}\ees
The existence and uniqueness of (\ref{3.1}) has been obtained by many authors, see
\cite{Cb,DZb}. Under the assumption the flux function $B$ is continuous with linear
growth. Debussche et al. \cite{DHV} obtained the following results, see Theorem 2.5 in
\cite{DMH}.

\begin{prop}\lbl{p3.1} There exists $((\tilde\Omega,\mathscr{\tilde F},\mathbb{\tilde P}),\tilde W,\tilde u)$ which
is a weak martingale solution to (\ref{3.1}) and for all $p\in[2,\infty)$ and $u_0\in L^p(\tilde\Omega;L^p)$,
    \bess
\tilde u\in L^p(\tilde\Omega;C([0,T];L^2);L^2)\cap L^p(\tilde\Omega;L^\infty(0,T;L^p))\cap
L^p(\tilde\Omega;L^2(0,T;W^{1,2})) .
   \eess
\end{prop}
Kim \cite{K2004} obtained the H\"{o}lder estimate of (\ref{3.1}), where
they used Bessel space similar to those in \cite{K1996,KpK,KK}. Based on the theory of
semigroup, Kuksin et al. \cite{Ku} obtained the H\"{o}lder estimate of (\ref{3.1}).

Let $D$ be an $A$-type bounded domain in $\mathbb{R}^{n+1}$.
Note that the Schauder estimate in this paper is interior estimate.
It is well known that the solution of
   \bess
u_t(t,x)=\Delta u+c(t,x)u+f(t,x)
    \eess
has the interior Schauder estimate if $c$ and $f$ are H\"{o}lder continuous.
Let $v$ be the solution of the following stochastic heat equation
    \bes
\left\{\begin{array}{lll}
du(t,x)=\Delta udt+g(t,x)dW(t),\ \ t>0,\,x\in\mathbb{R}^d,\\
u(0,x)=0,\ \  \,x\in\mathbb{R}^d.
  \end{array}\right.
  \lbl{3.3}\ees
Set $w:=u-v$, the $w$ satisfies that
      \bes
\left\{\begin{array}{lll}
w_t(t,x)=\Delta w+divB(u)+c(t,x)u+f(t,x),\ \ t>0,\,x\in\mathbb{R}^d,\\
u(0,x)=u_0(x),\ \  \,x\in\mathbb{R}^d,
  \end{array}\right.
  \lbl{3.4}\ees
Borrowing the idea from \cite{DMH} and using the \cite[Theorem 3.2]{DMH}, it is not hard
to prove that the solution $w$ of (\ref{3.4}) is H\"{o}lder continuous. That is,
there exists a positive constant $\gamma$ such that
    \bess
\mathbb{E}\|w\|_{C^\gamma(D_T)}=\mathbb{E}\sup_{t,x\in D_T}|u(t,x)|+
\mathbb{E}\sup_{(t,x)\neq(s,y)}\frac{|u(t,x)-u(s,y)|}{\max\{|t-s|^{\frac{1}{2}},|x-y|\}^\gamma}<\infty,
  \eess
where $D_T=[0,T]\times G$ and $G$ is a bounded domain in $\mathbb{R}^d$.
Note that
   \bess
\sup_{(t,x)\neq(s,y)}\frac{\mathbb{E}|u(t,x)-u(s,y)|}{\max\{|t-s|^{\frac{1}{2}},|x-y|\}^\gamma}\leq
\mathbb{E}\sup_{(t,x)\neq(s,y)}\frac{|u(t,x)-u(s,y)|}{\max\{|t-s|^{\frac{1}{2}},|x-y|\}^\gamma},
   \eess
we have the $w$ of (\ref{3.4}) belongs to $C^{\gamma}((\bar D_T;\delta);L^p(\Omega))$ for some $\gamma>0$.

It is easy to see that the mild solution $v$ of (\ref{3.3}) takes the following form
  \bess
v(t,x)=\mathcal {K}g(t,x)=\int_0^t\int_{\mathbb{R}^d}K(t,r,y)g(r,x-y)dydW(r),
  \eess
where $K(t,r;x,y)=(4\pi(t-r))^{-\frac{d}{2}}e^{-\frac{(x-y)^2}{4(t-r)}}$. It is easy to check that
the kernel function $K$ satisfies
   \bess
\int_{\mathbb{R}^d}K(t,r;x)dx=1,\ \ \int_{\mathbb{R}^d}|x|^\beta K(t,r;x)dx\leq N(T) \ \ for\ t\in[0,T],
   \eess
which implies that (\ref{2.3}) and (\ref{2.4}) with $\gamma_2=1$ hold.
Moreover, we have
  \bess
&&\int_0^s\left(\int_{\mathbb{R}^d}|K(t-r,z)-K(s-r,z)|(1+|z|^\beta)dz\right)^2dr\\
&=&(t-s)^2\int_0^s\left(\int_{\mathbb{R}^d}[\frac{d}{2(\xi-r)}-\frac{z^2}{4(\xi-r)^2}]
(4\pi(\xi-r))^{-\frac{d}{2}}e^{-\frac{z^2}{4(\xi-r)}}(1+|z|^\beta)dz\right)^2dr\\
&\leq& N(d,\beta)(t-s)^2\int_0^s(\xi-r)^{-d-2}\left(\int_{\mathbb{R}^d}[1+|z|^\beta+\frac{z^2}{4(\xi-r)}+\frac{|z|^{2+\beta}}{4(\xi-r)}]
e^{-\frac{z^2}{4(\xi-r)}}dz\right)^2dr\\
&\leq&N(d,T,\beta)(t-s)^2\int_0^s(\xi-r)^{-2}dr\\
&\leq&N(d,T,\beta)(t-s)^2[(\xi-s)^{-1}-\xi^{-1}]\\
&\leq&N(d,T,\beta,\theta)(t-s),
   \eess
where $\xi=\theta t+(1-\theta)s$ and $\theta\in(0,1)$. And thus (\ref{2.2}) holds with $\gamma_1=1$.
Therefore,
the assumptions of Theorems \ref{t2.1} and \ref{t2.2} hold.
It follows from Theorem \ref{t2.1} that
   \bess
v(t,x)\in C^{\beta}((\bar D_T;\delta);L^p(\Omega)).
  \eess
Combining the above results, we have the following
\begin{theo}\lbl{t3.1} Let $D_T$ be an $A$-type bounded domain in $\mathbb{R}^{d+1}$ such that $D_T\subset \mathcal {O}_T$.
Suppose the flux function $B$ is continuous with linear
growth, $u_0\in C^\beta(\mathbb{R}^d)$ and $g\in C^\beta(\mathbb{R}_+\times\mathbb{R}^d)$
with $g(0,0)=0$ almost surely, $0<\beta<1$, then the solution $u$ of (\ref{3.1}) is H\"{o}lder continuous in
domain $D_T$.
  \end{theo}
Similarly, we can use Theorem \ref{t2.2} to obtain the Schauder estimate of (\ref{3.1}), where $g$
does not depend on the time variable.

Next, we consider the following stochastic fractional heat equation
    \bes
\left\{\begin{array}{lll}
du(t,x)=\Delta^{\frac{\alpha}{2}} udt+g(t,x)dW(t),\ \ t>0,\,x\in\mathbb{R}^d,\\
u(0,x)=u_0(x),\ \  \,x\in\mathbb{R}^d,
  \end{array}\right.
  \lbl{3.5}\ees
where $\Delta^{\frac{\alpha}{2}}:=-(-\Delta)^{\frac{\alpha}{2}}$. Following the result of
\cite{XDLL}, the solution $u$ of (\ref{3.5}) can be written as
   \bes
u(t,x)&=&(\mathcal {G}\ast u_0)(t,x)+(\mathcal {G}\ast g)(t,x)\nm\\
&=&\int_{\mathbb{R}^d}p(t;x,y)u_0(y)dy+\int_0^t\int_{\mathbb{R}^d}p(t,r;x,y)g(r,y)dydW(r),
   \lbl{3.6}\ees
where the kernel function $p$ has the following properties:
 \begin{itemize}
 \item for any $t>0$,
\bess
\|p(t, \cdot)\|_{L^1(\R^d)}=1 \text{ for all } t>0.
\eess
 \item $ p(t, x, y)$ is  $C^\infty$ on $(0,\infty)\times \R^d\times \R^d$ for each $t>0$;

 \item   for $t>0$, $x, y\in \R^d$, $x\neq y$, the sharp estimate of $\widehat p(t, x)$ is
\bess
p(t, x, y)\approx \min\left( \frac{t}{|x-y|^{d+\alpha}}, t^{-d/\alpha}\right);
  \eess
\item  for $t>0$, $x, y\in \R^d$, $x\neq y$,   the  estimate of the first order derivative of $\widehat p(t, x)$ is
\bess
  |\nabla_x p(t, x, y)|\approx  |y-x|\min\left\{ \frac{t}{|y-x|^{d+2+\alpha}}, t^{-\frac{d+2}{\alpha}}\right\}.
  \eess
\end{itemize}

The notation $f(x)\approx g(x)$ means that there is a number $0<C<\infty$ independent of $x$,
i.e. a constant, such that for every $x$ we have $C^{-1}f(x)\leq g(x)\leq Cf(x)$.

\begin{prop}
\label{p3.2}{\rm\cite[Lemma 2.1]{XDLL}}
For any $m\geq 0$, we have
  \bess
\partial_x^m p(t, x)=\sum_{n=0}^{n=\lfloor \frac{m}{2}\rfloor}C_n |x|^{m-2n} \min \left\{ \frac{t}{|x|^{d+\alpha+2(m-n)}}, t^{-\frac{d+2(m-n)}{\alpha}}\right\},
     \eess
where $\lfloor \frac m 2\rfloor$ means the largest integer that is less  than $\frac m2$.
\end{prop}

By using Proposition \ref{p3.2}, we show the following

\begin{lem}\lbl{l3.1}Let $0\leq\epsilon<\frac{\alpha}{2}$. The following estimates hold.
   \bess
&&\int_0^s\left(\int_{\mathbb{R}^d}|\nabla^\epsilon p(t-r,z)-\nabla^\epsilon p(s-r,z)|(1+|z|^\beta)dz\right)^2dr\leq N(T,\beta)(t-s)^{\gamma},\\
&&\int_0^s\left(\int_{\mathbb{R}^d}|\nabla^\epsilon p(s-r,z)|dz\right)^2dr\leq N_0,\\
&& \int_s^t\left(\int_{\mathbb{R}^d}|\nabla^\epsilon p(t-r,z)|(1+|z|^\beta) dz\right)^2dr\leq N(T,\beta)(t-s)^{\gamma},
   \eess
where $\gamma=\frac{\alpha-2\epsilon}{\alpha}$.
\end{lem}

{\bf Proof.} For simplicity, we first prove the estimates with $\beta=0$ hold. It is
not hard to prove that when $\beta>0$, the index will be improved and the proof is omitted here.
Noting that $\partial_tp=-(-\Delta)^{\frac{\alpha}{2}}p:=\nabla^\alpha p$, when $\lfloor \frac {\alpha+\epsilon}{ 2}\rfloor<1$,
we have
    \bess
&&\int_0^s\left(\int_{\mathbb{R}^d}|\nabla^\epsilon p(t-r,z)-\nabla^\epsilon p(s-r,z)|dz\right)^2dr\\
&\leq&(t-s)^2\int_0^s\left(\int_{\mathbb{R}^d}|\nabla^{\alpha+\epsilon} p(\xi-r,z)|dz\right)^2dr\\
&\leq&(t-s)^2\int_0^s\left(\int_{\mathbb{R}^d}|z|^{\alpha+\epsilon}
 \min \left\{ \frac{\xi-r}{|z|^{d+3\alpha+2\epsilon}}, (\xi-r)^{-\frac{d+2\alpha+2\epsilon}{\alpha}}\right\}dz\right)^2dr\\
&\leq&(t-s)^2\int_0^s\left(\int_0^{(\xi-r)^{\frac{1}{\alpha}}}|z|^{\alpha+\epsilon}|z|^{d-1}(\xi-r)^{-\frac{d+2\alpha+2\epsilon}{\alpha}}d|z|\right.\\
&&\left.+\int_{(\xi-r)^{\frac{1}{\alpha}}}^\infty|z|^{\alpha+\epsilon}|z|^{d-1}\frac{\xi-r}{|z|^{d+3\alpha+2\epsilon}} d|z|\right)^2dr\\
&\leq&N(d,\alpha)(t-s)^2\int_0^s(\xi-r)^{-2\frac{\alpha+\epsilon}{\alpha}}dr\\
&\leq&N(d,\alpha,\theta)(t-s)^{\frac{\alpha-2\epsilon}{\alpha}},
   \eess
where $\xi=\theta t+(1-\theta)s$ and $\theta\in(0,1)$.

When $1\leq\lfloor \frac {\alpha+\epsilon}{ 2}\rfloor<2$, there is a little different from the
above discussion. Similarly, we get
   \bess
&&\int_0^s\left(\int_{\mathbb{R}^d}|\nabla^\epsilon p(t-r,z)-\nabla^\epsilon p(s-r,z)|dz\right)^2dr\\
&\leq&(t-s)^2\int_0^s\left(\int_{\mathbb{R}^d}|z|^{\alpha+\epsilon}
 \min \left\{ \frac{\xi-r}{|z|^{d+3\alpha+2\epsilon}}, (\xi-r)^{-\frac{d+2\alpha+2\epsilon}{\alpha}}\right\}dz\right)^2dr\\
&& +(t-s)^2\int_0^s\left(\int_{\mathbb{R}^d}|z|^{\alpha+\epsilon-2}
 \min \left\{ \frac{\xi-r}{|z|^{d+3\alpha+2\epsilon-2}}, (\xi-r)^{-\frac{d+2\alpha+2\epsilon-2}{\alpha}}\right\}dz\right)^2dr\\
&\leq&(t-s)^2\int_0^s\left(\int_0^{(\xi-r)^{\frac{1}{\alpha}}}|z|^{d-1}|z|^{\alpha+\epsilon-2}\left[|z|^{2}(\xi-r)^{-\frac{d+2\alpha+2\epsilon}{\alpha}}
+(\xi-r)^{-\frac{d+2\alpha+2\epsilon-2}{\alpha}}\right]d|z|\right.\\
&&\left.+\int_{(\xi-r)^{\frac{1}{\alpha}}}^\infty|z|^{\alpha+\epsilon-2}|z|^{d-1}\left[|z|^2\frac{\xi-r}{|z|^{d+3\alpha+2\epsilon}}
+\frac{\xi-r}{|z|^{d+3\alpha+2\epsilon-2}}\right] d|z|\right)^2dr\\
&\leq&N(d,\alpha)(t-s)^2\int_0^s(\xi-r)^{-2\frac{\alpha+\epsilon}{\alpha}}dr\\
&\leq&N(d,\alpha,\theta)(t-s)^{\frac{\alpha-2\epsilon}{\alpha}},
   \eess
where $\xi=\theta t+(1-\theta)s$ and $\theta\in(0,1)$.

Using Proposition \ref{p3.2} again, we have
   \bess
&&\int_0^s\left(\int_{\mathbb{R}^d}|\nabla^\epsilon p(s-r,z)|dz\right)^2dr\\
&\leq&\int_0^s\left(\int_{\mathbb{R}^d}|z|^\epsilon\min\left\{\frac{s-r}{|z|^{d+\alpha+2\epsilon}},(s-r)^{-\frac{d+2\epsilon}{\alpha}}\right\} dz\right)^2dr\\
&\leq&\int_0^s\left(\int_0^{(s-r)^{\frac{1}{\alpha}}}|z|^\epsilon(s-r)^{-\frac{d+2\epsilon}{\alpha}}|z|^{d-1}d|z|\right.\\
&&\left.\int_{(s-r)^{\frac{1}{\alpha}}}^\infty|z|^\epsilon\frac{s-r}{|z|^{d+\alpha+2\epsilon}} |z|^{d-1}d|z| \right)^2dr\\
&\leq&N(d)\int_0^s(s-r)^{-\frac{2\epsilon}{\alpha}}dr\\
&\leq& N(d,\alpha,\epsilon)s^{1-\frac{2\epsilon}{\alpha}}:=N_0<\infty\, .
   \eess

Similarly, we get
   \bess
&& \int_s^t\left(\int_{\mathbb{R}^d}|\nabla^\epsilon p(t-r,z)|(1+|z|^\beta) dz\right)^2dr\\
&\leq& N(d,\alpha,\epsilon)\int_s^t(t-r)^{-\frac{2\epsilon}{\alpha}} dr\\
&\leq&N(d,\alpha,\epsilon)(t-s)^{1-\frac{2\epsilon}{\alpha}}.
   \eess
The proof is complete. $\Box$

Theorem \ref{t2.1} implies that the solution $u$ of (\ref{3.6}) satisfying $u\in C^{\epsilon+\beta_1,\beta_1/2}((\bar D;\delta);L^p(\Omega))$,
where $\beta_1=\min\{\beta,2\gamma\}$.
\begin{theo}\lbl{t3.2} Let $D_T$ be a $A$-type bounded domain in $\mathbb{R}^{d+1}$ such that $D_T\subset \mathcal {O}_T$.
Suppose that $u_0\in C^\beta(\mathbb{R}^d)$ and $g\in C^\beta(\mathbb{R}_+\times\mathbb{R}^d)$
with $g(0,0)=0$ almost surely, $0<\beta<1$, then the solution $u$ of (\ref{3.5}) is H\"{o}lder continuous in
domain $D_T$.
  \end{theo}

\begin{remark}\lbl{r3.1} Comparing with Theorems \ref{t3.1} and \ref{t3.2}, we find that if we take $\epsilon=0$, then
Theorem \ref{t3.2} with $\alpha=2$ becomes Theorem \ref{t3.1}. Let we compare the index of spatial variable. Theorem
\ref{t3.1} shows that the index is $\beta$ and Theorem
\ref{t3.2} shows that the index is $\epsilon+\min\{\beta,2\gamma\}$. When $\beta\leq2\gamma$, the result of
Theorem \ref{t3.2} is better than that of Theorem \ref{t3.2}.
\end{remark}

\subsection{Application to Fractional Heat Equations Driven by L\' evy Noise}

For simplicity, we consider the following SPDEs
    \bes
\left\{\begin{array}{lll}
du(t,x)=\Delta^{\frac{\alpha}{2}} u(t,x)dt+\ds\int_Zg(t,x,z)\tilde N(dt,dz),\ \ t>0,\,x\in\mathbb{R}^d,\\
u(0,x)=u_0(x),\ \  \,x\in\mathbb{R}^d,
  \end{array}\right.
  \lbl{3.7}\ees
where $\Delta^{\frac{\alpha}{2}}=-(-\Delta)^{\frac{\alpha}{2}}$. The well-posedness of
(\ref{3.7}) has been proved by \cite{KK}. The solution of (\ref{3.7}) can be written as
    \bes
u(t,x)&=&(\mathcal {G}\ast u_0)(t,x)+(\mathcal {G}\ast g)(t,x)\nm\\
&=&\int_{\mathbb{R}^d}p(t;x,y)u_0(y)dy+\int_0^t\int_{\mathbb{R}^d}\int_Zp(t,r;x,y)g(r,y,z) dy\tilde N(dt,dz).
   \lbl{3.8}\ees
Using the properties of $g$ and Lemma \ref{l3.1}, it is
easy to verify that all the assumptions in Theorem \ref{2.3} hold for the kernel function.
\begin{theo}\lbl{t3.3} Suppose that $u_0\in C^\beta(\mathbb{R}^d)$ with $\beta<\alpha$ and the function $g$ satisfies that
    \bess
  |g(t,x,z)-g(s,y,z)|\leq C_g\max\left\{(t-s)^{\frac{1}{2}},|x-y|\right\}^\beta g_1(z),\ \ \
  {\rm for \ all}\ \ z\in Z,\ \ a.s.,
      \eess
and $g(0,0,z)=0$ uniformly for $z\in Z$ almost surely, where there exists a constant $p_0>1$ such that $g_1(z)$ satisfies that
    \bess
\mathbb{E}\left[\left(\int_Z|g_1(z)|^2\nu(dz)\right)^{p_0/2}+\int_Z|g_1(z)|^{p_0}\nu(dz)\right]<\infty.
  \eess
Let $D$ be a $A$-type bounded domain in $\mathbb{R}^{d+1}$ such that $\bar D\subset \mathcal {O}_T$.
Then the solution $u$ of (\ref{3.7}) is H\"{o}lder continuous in
domain $D_T$.
  \end{theo}

\medskip

\noindent {\bf Acknowledgment} The first author was supported in part
by NSFC of China grants 11301146, 11531006, 11501577, 11171064.

\begin{appendix}
\section{Appendix}
\begin{lem}
A Companato space defined by Definition 2.1 is a Banach space. Moreover, if $1\leq p\leq q<\infty$, $(\theta-p)/p\leq(\sigma-p)/q$, it holds that
\bess
\mathscr{L}^{q,\sigma}(D;\delta)\subset\mathscr{L}^{p,\theta}(D;\delta).
\eess
\end{lem}
{\bf Proof.} Let $\{u^n\}_{n\geq 1}$ be a Cauchy sequence in $\mathscr{L}^{p,\theta}(D;\delta)$, then $\{u^n\}_{n\geq 1}$ is also a Cauchy in $L^p(D)$. Therefore, there is a measurable function $u\in L^p(D)$, such that
\bess
u^n\longrightarrow u \ \  {\rm in} \ \ L^p(D) .
\eess
By virtue of the Fato lemma, we have the following estimate:
\bess
[u-u^n]_{\mathscr{L}^{p,\theta}(D;\delta)}^p&=&\sup_{X\in D,d\geq\rho>0}\frac{1}{|D(X,\rho)|^\theta}
\int_{D(X,\rho)}|(u-u^n)(Y)-(u-u^n)_{X,\rho}|^pdY
\cr
&=&\sup_{X\in D,d\geq\rho>0}\lim_m\frac{1}{|D(X,\rho)|^\theta}
\int_{D(X,\rho)}|(u^m-u^n)(Y)-(u^m-u^n)_{X,\rho}|^pdY
\cr&\leq& \limsup_m\sup_{X\in D,d\geq\rho>0}\frac{1}{|D(X,\rho)|^\theta}
\int_{D(X,\rho)}|(u^m-u^n)(Y)-(u^m-u^n)_{X,\rho}|^pdY.
\eess
To prove $u^n\rightarrow u$ in $\mathscr{L}^{p,\theta}(D;\delta)$, it suffices to show $[u]_{\mathscr{L}^{p,\theta}(D;\delta)}<\infty$. And this estimate holds clearly by using above calculation, so the first part of the lemma is finished.

To verify the second part, we use the H\"{o}lder inequality for $1\leq p\leq q<\infty$, $0<\theta,\sigma$ to gain that: for every $u\in \mathscr{L}^{q,\sigma}(D;\delta)$
\bess
[u]_{\mathscr{L}^{p,\theta}(D;\delta)}&=&\left(\sup_{X\in D,d\geq\rho>0}\frac{1}{|D(X,\rho)|^\theta}
\int_{D(X,\rho)}|u(Y)-u_{X,\rho}|^pdY\right)^{1/p}
\cr&\leq&\left(\sup_{X\in D,d\geq\rho>0}\frac{1}{|D(X,\rho)|^{\theta q/p+p-q}}
\int_{D(X,\rho)}|u(Y)-u_{X,\rho}|^qdY\right)^{1/q}
\cr&=&\left(\sup_{X\in D,d\geq\rho>0}\frac{1}{|D(X,\rho)|^{\sigma}}\frac{1}{|D(X,\rho)|^{\theta q/p+p-q-\sigma}}
\int_{D(X,\rho)}|u(Y)-u_{X,\rho}|^qdY\right)^{1/q}.
\eess
If $\theta q/p+p-q-\sigma\leq 0$, i.e. $(\theta-p)/p\leq(\sigma-p)/q$, then
\bess
[u]_{\mathscr{L}^{p,\theta}(D;\delta)}\leq C\left(\sup_{X\in D,d\geq\rho>0}\frac{1}{|D(X,\rho)|^{\sigma}}
\int_{D(X,\rho)}|u(Y)-u_{X,\rho}|^qdY\right)^{1/q},
\eess
for $D$ is bounded.

On the other, $L^q(D)\subset L^p(D)$, we complete the proof. $\Box$

\end{appendix}
 \end{document}